\documentclass{article}
\usepackage{amsfonts,amssymb,a4wide,amsmath}
\newcommand{\R}{\mathbb{R}}
\newcommand{\N}{\mathbb{N}}
\newcommand{\C}{\mathbb{C}}
\newcommand{\Z}{\mathbb{Z}}
\newtheorem{prop}{Proposition}[section]
\newtheorem{teo}[prop]{Theorem}
\newenvironment{pro}{\noindent{\bf Proof. }}{\qed}
\newcommand{\qed}{\nobreak \hfill \rule{.20cm}{.20cm} \vspace*{.5cm}}
\begin{document}

\begin{center}

{\bf  IMAGE NORMALIZATION OF WIENER-HOPF OPERATORS AND\\
BOUNDARY-TRANSMISSION VALUE PROBLEMS\\
FOR A JUNCTION OF TWO HALF-PLANES}

\vspace{5mm}

{A. MOURA SANTOS$^*$}, {N.J. BERNARDINO $^{**}$}

$^*$\it{Department of Mathematics, Instituto Superior T\'ecnico, Technical University of Lisbon,\\}
\it{1049-001 Lisbon, Portugal\\}
E-mail: amoura@math.ist.utl.pt

$^{**}$\it{ESB 2440-062 Batalha, Portugal\\}

\end{center}

\begin{abstract}
The present paper deals with an application  of the image normalization technique for certain classes of Wiener-Hopf operators (WHOs) associated  to ill-posed boundary-transmission value problems. We briefly describe the method of normalization and then apply it to boundary-transmission value problems issued from diffraction problems for a junction of two half-planes, which are relevant in mathematical physics applications. We consider different boundary-transmission conditions on the  junction of the two semi-infinite half-planes and analyze the not normally solvability of the corresponding operators.\end{abstract}

{\bf keywords}: Diffraction by a junction of two half-planes; Boundary-transmission value problems; Not normally solvable operators; Image normalization; Wiener-Hopf operators.

%\bodymatter

\section{Introduction}

We are interested in operators, namely Wiener-Hopf operators (WHOs), which arise in the context of diffraction problems of electromagnetic and acoustic waves and are strongly related to the operator description of the corresponding boundary-transmission value problems. In general, for many relevant physical situations of the boundary the corresponding boundary value problems are ill-posed, i.e. the associated operators are not normally solvable, see e.g. \cite{MeisterSpeck89, SantosTeixeira95, AmSpeck97, AmSpeckTeixeira98}. This was the main reason why one of the authors of the present paper developed in her PhD work a method of image normalization in order to convert not normally solvable WHOs into operators with closed image. This method was firstly applied to boundary value problems on the half-plane  \cite{AmSpeckTeixeira98}, but can be successfully used for other geometries of the boundary. In this paper we describe how to apply it for a junction of two half-planes, for a strip see e.g. \cite{AmCastro04}.  The method of image normalization is one of the possible ways of normalizing  bounded linear operators acting between Banach spaces \cite{Kravchenko85}, and works very efficiently for the operators under consideration. 

In Section \ref{aba:sec1} we describe the class of  boundary-transmission value problems which arise from the diffraction of a plane wave by a junction of two half-planes. Diffraction by a two-part plane is relevant for many practical applications, see e.g. \cite{Norris95, Rojas89,  Serbest95}.   Starting from the standard operator procedure of the classical survey of Meister and Speck \cite{MeisterSpeck89}, we associate with the physical problem an operator and then  prove the equivalence of this operator to a WHO. It is also in Section \ref{aba:sec1} that we introduce same needed notation. Section \ref{image} is dedicated to summarize the method of image normalization of the WHOs under consideration. There we present the main results of   \cite{AmSpeckTeixeira98} without proofs and using a notation more convenient for our present purposes. The next three  sections describe and analyze chosen examples of image normalization of WHOs coming from different boundary-transmission conditions on the two half-planes. For instance, in Section \ref{example1} we first derive the WHO for boundary-transmission conditions of arbitrary orders on the two banks of the two half-planes, and then consider the image normalization when all four orders are even. Section \ref{example2} is devoted to consider  boundary-transmission conditions with only normal derivatives of the same order on the upper and lower banks of the two half-planes, respectively. Finally, in Section \ref{example3} we consider a simpler boundary condition on the left half-plane and a boundary-transmission condition with oblique derivatives on the right half-plane.

\section{Boundary-transmission problems and WHOs}\label{aba:sec1}

In order to study the WHOs, we begin with the formulation of the following general boundary-transmission value problem, we call it the Problem ${\mathcal P}$, for the diffraction of a plane wave by a junction of two half-planes in the natural setting of locally finite energy norm.

{\bf Problem ${\mathcal P}$}. Find 
 ${\varphi} \in L^2({\R}^2)$, with
${\varphi}_{|{\R\times \R_{\pm}}}= {\varphi}^{\pm}\in H^{1}({\R\times \R}_{\pm})$, s.t.
\begin{equation}  \label{Helmholtz}
 \begin{array}{ccl}
\left(\Delta + k_0^{2}\right) {\varphi}^{\pm} \;\; = \;\; 0 \quad &\mbox{in}&\quad
{{\R\times \R}}_{\pm}\,,
 \end{array}
\end{equation}
\begin{equation}  \label{boundary1}
 B_j^- \varphi(x) =\sum_{\sigma_{1} +\sigma_{2}\leq m_j}
  a_{\sigma, j}^+ (D^{\sigma}\varphi^+ )(x,0)+a_{\sigma, j}^-
(D^{\sigma}\varphi^-)(x,0) = h_j(x)\,\,\, {\rm on} \,\,\, \R_-,  
\end{equation}
\begin{equation}   \label{boundary2}
 B_j^+ \varphi(x) =\sum_{\sigma_{1} +\sigma_{2}\leq m'_j}
  b_{\sigma, j}^+ (D^{\sigma}\varphi^+ )(x,0)+b_{\sigma, j}^-
(D^{\sigma}\varphi^-)(x,0) = g_j(x)\,\,\, {\rm on} \,\,\, \R_+,     
\end{equation}
where $\sigma=(\sigma_1,\sigma_2)$, $\sigma_j \in \N_0$, and $m=(m_1,m_2)$, $m'= (m'_1,m'_2)$, $m_j, m'_j \in \N_0$, represent the order of the boundary operators $B_j^{-}$ and $B_j^{+}$, respectively,  with $j=1, 2$ corresponding to the upper and lower banks of  both left $\R_-$ and right $\R_+$ half-lines  \footnote{As a consequence of the physics of the wave diffraction the boundaries, i.e. the two half-planes, can be identified with these two subsets of the real line.}. The coefficients  $a_{\sigma, j}^{\pm}, b_{\sigma, j}^{\pm} \in \C$ simulate physical properties of the boundaries. 
For instance,  
for $\sigma=(0,1)$, $m_j=(0,1)$, $a_{\sigma, j}^{+}=1, a_{\sigma, j}^{-}=-1$ and $h=(h_1,h_2)=(0,0)$ in (\ref{boundary1}), $B_j^-$ consists of the trivial jump of Dirichlet data and Neuman data on $\R_-$ and is usually known as the transmission condition for the Sommerfeld problem, see e.g. \cite{MeisterSpeck89}. On the other hand for  $\sigma_{1} +\sigma_{2}\leq 1$ in (\ref{boundary2}), $B_j^+$ consists of a linear combination of Dirichlet, Neumann, and oblique derivative data as considered in \cite{AmSpeck97}. It is also physically meaningful to consider linear combinations of higher derivatives \cite{SantosTeixeira95} both normal and tangential, with coefficients  $a_{\sigma, j}^{\pm}, b_{\sigma, j}^{\pm} \in \C$ depending on the materials of the boundary. In the Helmholtz equation (\ref{Helmholtz}) $k_0$ stands for the complex wave number with positive real and imaginary part, i.e. ${\rm Re}\,k_0>0$ and ${\rm Im}\,k_0>0$.

Following the operator procedure of the classical survey of Meister and Speck \cite{MeisterSpeck89}, we describe the Problem ${\mathcal P}$ by a unique equation 
\begin{equation}\label{operatorP}{\mathcal P}\varphi=g,
\end{equation}
where ${\mathcal P}$ is a linear operator {\it associated} to the Problem  ${\mathcal P}$ which acts like
$${\mathcal P}: D({\mathcal P}) \rightarrow H^{1/2-m'_1} (\R_+)\times H^{1/2-m'_2} (\R_+).$$
The domain $D({\mathcal P})$ is given by the subspace of $H^{1}({\R\times \R}_{\pm})$ whose elements fulfill the Helmholtz equation (\ref{Helmholtz}) and  the boundary-transmission conditions in (\ref{boundary1}) and the image space is characterized by $g=(g_1,g_2) \in  H^{1/2-m'_1} (\R_+)\times H^{1/2-m'_2} (\R_+)$ according to the trace theorem (see below) and the representation formula applied to (\ref{boundary2}) with $m'=(m'_1,m'_2)$.  

Then the next goal is to prove the equivalence relation ${\mathcal P}=EWF$, where $E$ and $F$ are bounded invertible linear operators and the operator $W$ acts between ${H}^r_+=H^{r_1}_+\times H^{r_2}_+$ subspaces of $H^r$ distributions supported on $\overline{\R}_+$, and
$H^s(\R_+)=H^{s_1}(\R_+)\times H^{s_2}(\R_+)$ restrictions of $H^s$  distributions on ${\R_+}$. Each component of $H^s$ (or $H^r$) is a Bessel potential space of order $s_j$ (or $r_j$) defined by
$$H^{s_j}=\{\phi \in {\mathcal S}':  {\mathcal F}^{-1}
(\xi^2+1)^{s_j/2}\cdot {\mathcal F}\phi \in L^2\}$$
where $ {\mathcal F}$ represents the Fourier transformation \footnote{In the context of Problem ${\mathcal P}$, it is  also commonly used  $(\xi^2-k_0^2)^{s_j/2}$ in the definition of the Bessel potential spaces, in which case the branch cuts are defined along $\pm k_0\pm i\epsilon$, $\epsilon\geq 0$.}. These are well-known Hilbert spaces.The topologies are the usually subspace topology for each component $H^{r_j}_+$ and the quotient space topology for each $H^{s_j}(\R_+)$.

Now we describe in more detail how to obtain from (\ref{Helmholtz})-(\ref{boundary2}) the WHO equivalent to the {\it associated} operator ${\mathcal P}$.

We start with the standard representation formula, see e.g. \cite{Speck86}, for the solutions of the Helmholtz equation (\ref{Helmholtz})
\begin{equation} \label{represent}
\varphi(x,y)={\mathcal K}\varphi_0(x,y)={\mathcal F}^{-1}_{\xi \rightarrow x}\left\{e^{-\beta(\xi)y}\hat{\varphi}_0^+(\xi)\chi_{+}(y)+e^{\beta(\xi)y}\hat{\varphi}_0^-(\xi)\chi_{-}(y)\right\},
\end{equation} 
where $\varphi_0=({\varphi}_0^+,{\varphi}_0^- ) \in H^{1/2}\times H^{1/2}$ is the trace vector of $\varphi^{\pm}$ due to the banks of $\R_{\pm}$, $\hat{\varphi}_0^{\pm}$ represent the Fourier transform of the traces, $\chi_{\pm}$ denote the characteristic function of the positive and negative half-line, respectively, and $\beta(\xi)=\sqrt{\xi^2-k_0^2}$. For the operator ${\mathcal K}$ in (\ref{represent})  the following result holds.

\begin{teo} \label{traceoperator}
Let $B_{-}= {\mathcal F}^{-1}
{\Phi_{-}} \cdot {\mathcal F}$ be the linear bounded operator
$$
B_{-}:[H^{1/2}]^2 \rightarrow  H^{1/2-m_1}\times H^{1/2-m_2}
$$
with Fourier symbol
\begin{equation} \label{Fourier1}
\Phi_-=
\left[
   \begin{array}{cc}
 {\underset{ |\sigma|\leq m_{1}}\sum}
          a_{\sigma, 1}^+ (-i\xi)^{\sigma_1}(-\beta(\xi))^{\sigma_2}
  & 
 {\underset{ |\sigma|\leq m_{1}}\sum}
a_{\sigma, 1}^- (-i\xi)^{\sigma_1}\beta(\xi)^{\sigma_2}   \\  
 {\underset{ |\sigma|\leq m_{2}}\sum}  a_{\sigma, 2}^+ (-i\xi)^{\sigma_1}(-\beta(\xi))^{\sigma_2} & 
   {\underset{ |\sigma|\leq m_{2}}\sum}a_{\sigma, 2}^- (-i\xi)^{\sigma_1}\beta(\xi)^{\sigma_2}  \end{array}
 \right]  
\end{equation}
s.t.  $\rm {det}\Phi_-\neq 0$. Furthermore, consider the restricted operator $\tilde{B}_-={\rm {Rst}} B_-: [H^{1/2}]^2 \rightarrow  H_+^{1/2-m_1}\times H_+^{1/2-m_2}$ with Fourier symbol also given by (\ref{Fourier1}), but acting into the $H_+^s$ spaces.

Then the operator ${\mathcal K}$ in (\ref{represent}) is invertible by the trace operator  $T_0:D({\mathcal P}) \rightarrow Y_0 $, where $D({\mathcal P})$ is a closed subspace of the direct sum $H^1(\R\times \R_+)\bigoplus H^1(\R\times \R_-)$ and the image space is given by
\begin{equation} \label{imagetrace}
Y_0=\left\{  \varphi_0=(\varphi_0^+,\varphi_0^-) \in [H^{1/2}]^2:    {\mathcal F}^{-1}
{\Phi_{-}} \cdot {\mathcal F}\varphi_0 -\ell ^{(c)} h \in  H_+^{1/2-m_1}\times H_+^{1/2-m_2}  \right\}
\end{equation}
where $h=(h_1,h_2)$ is the data from (\ref{boundary1}), $\ell ^{(c)}$ represents a continuous extension operator, of even type for $m_j$ even, and odd type for $m_j$ odd, and left invertible by the restriction operator $r_{+}$. 

Moreover, for $\varphi_0=({\varphi}_0^+,{\varphi}_0^- )= \tilde{B}_-^{-1}\left(v^{+}+\ell ^{(c)} h\right)$ with $v^{+}\in  H_+^{1/2-m_1}\times H_+^{1/2-m_2}$ in (\ref{represent}), 
the operator $\tilde{B}_-T_0$ is continuously invertible by ${\mathcal K}\tilde{B}_-^{-1}$. 
\end{teo}

\begin{pro}
The trace operator $T_0:D({\mathcal P}) \rightarrow Y_0$ is here defined as an operator that acts between spaces of order greater or equal one half. For these orders of the spaces we have automatically surjectivity and right invertibility. Then the left invertibility is obtained by choosing the space $Y_0$ in (\ref{imagetrace}) as a subspace of order one half and defined in such a way that it contains zero extensions of the corresponding trace values that appear in our problem, see e.g. \cite{AmCastro04} for a discussion when this fails. Therefore, we have invertibility of $T_0$ by the operator ${\mathcal K}$ given by (\ref{represent}). 

Consider now $\varphi_0= \tilde{B}_-^{-1}\left(v^{+}+\ell ^{(c)} h\right)$ in the representation formula (\ref{represent}). Then $\tilde{B}_-\varphi_0=v^{+}+\ell ^{(c)} h$ and we have
$$
\tilde{B}_-T_0{\mathcal K}\tilde{B}_-^{-1}(v^{+}+\ell ^{(c)} h)=\tilde{B}_-T_0{\mathcal K}\varphi_0=\tilde{B}_-\varphi_0=v^{+}+\ell ^{(c)} h,
$$
and also
$$
{\mathcal K}\tilde{B}_-^{-1}\tilde{B}_-T_0\varphi={\mathcal K}\tilde{B}_-^{-1} \tilde{B}_-\tilde{B}_-^{-1}(v^{+}+\ell ^{(c)} h)={\mathcal K}\left(\tilde{B}_-^{-1}(v^{+}+\ell ^{(c)} h)\right)={\mathcal K}\varphi_0=\varphi,$$ 
i.e.  $\tilde{B}_-T_0$ is continuously invertible by ${\mathcal K}\tilde{B}_-^{-1}$
\end{pro}

Let us now prove the following equivalence result.

\begin{teo} \label{equivalencetoWHO}
Let $B_{+}= {\mathcal F}^{-1}
{\Phi_{+}} \cdot {\mathcal F}$ be the linear bounded operator
$$
B_{+}: [H^{1/2}]^2 \rightarrow  H^{1/2-m'_1}\times H^{1/2-m'_2}
$$
with Fourier symbol
\begin{equation} \label{Fourier2}
\Phi_+=
\left[
   \begin{array}{cc}
 {\underset{ |\sigma|\leq m'_{1}}\sum}
          b_{\sigma, 1}^+ (-i\xi)^{\sigma_1}(-\beta(\xi))^{\sigma_2}
  & 
 {\underset{ |\sigma|\leq m'_{1}}\sum}
b_{\sigma, 1}^- (-i\xi)^{\sigma_1}\beta(\xi)^{\sigma_2}   \\  
 {\underset{ |\sigma|\leq m'_{2}}\sum}  b_{\sigma, 2}^+ (-i\xi)^{\sigma_1}(-\beta(\xi))^{\sigma_2} & 
   {\underset{ |\sigma|\leq m'_{2}}\sum}b_{\sigma, 2}^- (-i\xi)^{\sigma_1}\beta(\xi)^{\sigma_2}  \end{array}
 \right]  
\end{equation}
s.t.  $\rm {det}\Phi_+\neq 0$. Moreover, let the conditions of Theorem \ref{traceoperator} hold. Then, the operator ${\mathcal P}$ in (\ref{operatorP}) is equivalent to the WHO
\begin{equation} \label{equivalentWHO}
W= r_{+}B_+\tilde{B}_-^{-1}  : H^{1/2-m_1}_+ \times H^{1/2-m_2}_+ \rightarrow H^{1/2-m'_1}(\R_+) \times H^{1/2-m'_2}(\R_+),
\end{equation} 
with Fourier symbol $\Phi=\Phi_+\Phi_-^{-1}$. The equivalence relation is given by 
\begin{equation} \label{equivalentrelation}
{\mathcal P}=W\tilde{B}_- T_0,
\end{equation}
i.e. the operators ${\mathcal P}$ and $W$ coincide up to bijective factors.
\end{teo}

\begin{pro}
Since from  Theorem \ref{traceoperator}, we have $\varphi_0=\tilde{B}_-^{-1}\left(v^{+}+\ell ^{(c)} h\right)$, it follows that $v^{+}= \tilde{B}_-\varphi_0-\ell ^{(c)} h$ and
$$
Wv^{+}=W( \tilde{B}_-\varphi_0 -\ell ^{(c)} h)= r_{+}B_+\tilde{B}_-^{-1} (\tilde{B}_-\varphi_0 -\ell ^{(c)} h)=r_+B_+\varphi_0-r_+B_+\tilde{B}_-^{-1} \ell ^{(c)} h=g-r_+B_+\tilde{B}_-^{-1} \ell ^{(c)} h.
$$
On the other hand, assuming that (\ref{equivalentrelation}) holds, we can write
$$
{\mathcal P}\varphi=W\tilde{B}_- T_0\varphi=W\tilde{B}_-\varphi_0= W\tilde{B}_-\tilde{B}_-^{-1} (v^+ +\ell ^{(c)} h)=Wv^++W\ell ^{(c)} h,$$
which after substituting $Wv^+$ by the expression obtained before, gives ${\mathcal P}\varphi=g$. This proves the equivalence between the two operators with the equivalence relation given by (\ref{equivalentrelation}).
\end{pro}

Let us study in more detail the general structure of the obtained operators. First, we formally rewrite the WHO as an operator
\begin{equation} \label{matrixWHO}
W= r_{+}A\mid_{H^r_+}  : H^r_+ \rightarrow H^s(\R_+), 
\end{equation}
with $A= {\mathcal F}^{-1}
{\Phi}\cdot {\mathcal F}$, $\Phi=\Phi_+\Phi_-^{-1}$, constituting a translation invariant homeomorphism
with a matrix  {\it Fourier symbol} $\Phi \in L_{\rm loc}^{\infty}$. Note that the elements of the Fourier symbol $\Phi$ of $W$ in (\ref{equivalentWHO}), due to (\ref{Fourier1}) and (\ref{Fourier2}), and given arbitrary orders $m_j,m'_j$, and coefficients $a_{\sigma, j}^{\pm},b_{\sigma, j}^{\pm}$ are rational functions of $\xi$ and $\beta(\xi)=\sqrt{\xi^2-k_0^2}$, see next sections for details. Then, lifting the WHO $W$ into $L^2$ see e.g. \cite{Eskin81}, we obtain the {\it lifted} WHO

\begin{equation}  \label{WHO_0}
W_0= r_{+}A_0\mid_{[L^2_+]^n}  : [L^2_+]^2 \rightarrow [L^2(\R_+)]^2,
\end{equation}
where 
$A_0= {\mathcal F}^{-1}
{\Phi}_0\cdot {\mathcal F}$, $\Phi_0 \in L^{\infty}(\R)^{2\times 2}$. In this paper we assume first that $\Phi_0 \in \mathcal{G}C^{\nu}(\ddot{\R})^{2\times 2}$  i.e. that the lifted Fourier symbol belongs to the invertible algebra of H\"{o}lder continuous $2\times 2$ matrix functions defined on $\ddot{\R}=[-\infty,+\infty]$. In fact the elements of the lifted Fourier symbol $\Phi_0$ are bounded rational functions of  $\rho(\xi)=\sqrt{\frac{\xi-k_0}{\xi+k_0}}$ and of $\xi \beta(\xi)^{-1}$, and we should always assume first ${\rm det}\, \Phi_0(\xi)\neq 0$, $\xi \in \ddot{\R}$, in order to get normal type WHOs, and then concentrate on the study of not normally solvability  \footnote{In this paper we are interested in the not normally solvable cases, thus we assume first that the coefficients in (\ref{boundary1})-(\ref{boundary2}) are s.t. ${\rm det}\Phi_0\neq 0$, i.e. the matrix does not degenerate on $\ddot{\R}$, and for these coefficients analyze the case of  not normally solvable WHOs.}.

The following Fredholm criterium is well-known \cite{MikhlinProssdorf86} for the lifted WHO in (\ref{WHO_0}). 
The operator $W_0$ is normally solvable iff
\begin{equation} \label{criterium}
{\rm det}(\mu \Phi_0(-\infty)+(1-\mu) \Phi_0(+\infty))\neq 0\,,\,\,\,\mu \in ]0,1[.
\end{equation}

As a consequence of the limits $\rho(\xi)\rightarrow \pm 1$ and $\xi\beta(\xi)^{-1}\rightarrow \pm 1$ as $\xi\rightarrow \pm \infty$, this condition does not hold for a large class of WHOs in (\ref{equivalentWHO}), and from the the equivalence relation (\ref{equivalentrelation}), the same is true for the associated  operator ${\mathcal P}$. Therefore there is a need to achieve the image normalization of both operators.

Finally, we shall also use the zero extension operator $\ell^{(0)}$ and the following Bessel potential operators \cite{AmSpeckTeixeira98} for $w \in \C, k_0 \in \C$, ${\rm Im}\,k_0>0$

$$\Lambda_{\pm}^w= {\mathcal F}^{-1}
{\lambda_{\pm}^w}\cdot {\mathcal F}:H^s \rightarrow H^{s-{\rm Re} w},$$ 
where we introduced $\lambda_{\pm}(\xi)=\xi\pm k_0$, a notation often used in this context.
%%%%%%%%%%%%%%%%%%%%%%%%%

\section{Image normalization of WHOs in scalar and matrix cases} \label{image}

We briefly describe the main results of our approach (for proofs see  \cite{AmSpeckTeixeira98}) towards the normalization of the WHOs defined by (\ref{matrixWHO}) with the corresponding lifted Fourier symbol $\Phi_0 \in \mathcal{G}C^{\nu}(\ddot{\R})^{2\times 2}$ for which the Fredholm criterium (\ref{criterium}) doesn't hold.  The method is based on two central ideas: first, we want the domain of the operator to remain a space of locally finite energy, and second, we change the image space in a {\it minimal} way. The following scalar result  \cite{AmSpeckTeixeira98} helps to understand the method for the matrix case.

\begin{teo} \label{Main Theorem-scalar}
Let us consider the scalar WHO of normal type, which acts symmetrically, i.e. $r=s$ 

$$W_s=W_s(\Phi)= r_{+}A\mid_{H^s_+}  : H^s_+ \rightarrow H^s(\R_+).$$
Then for the {\it critical} orders \cite{Duduchava79} $s+\eta+1/2\in \Z$, where $\eta=\frac{1}{2\pi i}\int_{\R} d \arg \Phi$, the operator $W_s$ is  not normally solvable.

Introducing $w=\eta+i \tau$,  with $\tau=\frac{1}{2\pi}\ln |\Phi(-\infty)/\Phi(+\infty)|$, we define
the {\it image normalized} operator $\breve{W}_s$ by

$$\breve{W}_s={\rm Rst}W_s: H^s_+ \rightarrow \breve{H}^{s-i\tau}(\R_+)$$
where  $\breve{H}^{s-i\tau}(\R_+)=r_+\Lambda_-^{-s+i\tau-1/2}H_+^{-1/2}\subset H^{{\rm Re}w}(\R_+)$.The image space of $\breve{W}_s$ solves the {\it normalization problem} for $\{W_s=W_s(\Phi):\Phi \in \mathcal{G}C^{\nu}(\ddot{\R})\,,\,\, \nu \in ]0,1[\,,\,\,{\rm im}W_s\neq \overline{{\rm im}W_s}\}$.
\end{teo}

The normalization in the matrix case is based on the same idea of using the jump at infinity of the lifted Fourier symbol to change the image space in a {\it minimal} way. The following result can be found in  \cite{AmSpeckTeixeira98} for the $n\times n$ matrix case, but here we state it for the $2\times 2$ matrix case, which will be enough for our purposes.

\begin{teo} \label{Main Theorem-matrix}
Consider the WHO $W$ in (\ref{matrixWHO}),
with $r=(r_1,r_2)$, $s=(s_1,s_2)$ and s.t. the corresponding lifted Fourier symbol $\Phi_0 \in \mathcal{G}C^{\nu}(\ddot{\R})^{2\times 2}$. Moreover, let this Fourier symbol $\Phi_0$ has a jump at infinity. Assume that $\lambda_1,\lambda_2$ are the eigenvalues of $\Phi_0$ and write

\begin{equation} \label{jump}
\Phi_0^{-1}(-\infty)\Phi_0(+\infty)=T^{-1}{\rm diag}(\lambda_1,\lambda_2)T,
\end{equation}
where $T\in \mathcal{G}\C^{2\times 2}$, and moreover assume that $\lambda_1=e^{2i\pi w_1}$ with ${\rm Re}w_1=-1/2$ ($\lambda_2=e^{2i\pi w_2}$ with ${\rm Re}w_2\neq-1/2$). Then the WHO $W$ is not normally solvable for the given $s_1$ and $\eta_1=-1/2$, but we can define the {\it image normalized} operator $\breve{W}$ by

$$\breve{W}: H^r_+ \rightarrow Y_1=r_+\Lambda_-^{-s}T\ell^{(0)}\{  \breve{H}^{-i\tau}(\R_+)\times L^2(\R_+)\}$$
with $\tau=\frac{1}{2\pi}\ln |\Phi_0(-\infty)/\Phi_0(+\infty)|$ which corresponds to the eigenvalue $\lambda_1$. The image space $Y_1$ of the restricted WHO $\breve{W}$ solves the {\it normalization problem} for $W$.
\end{teo}

Remark that we say that $Y_1$ solves the normalization problem for the WHO  and denote by $\breve{W}$ the corresponding {\it image normalized} operator.

In Theorem \ref{Main Theorem-matrix} we assumed that the eigenvalues  of $\Phi_0$ are different and that only one of the eigenvalues, $\lambda_1$ in (\ref{jump}), is responsible for the jump at infinity. As we will see in next sections, it is possible that we get an eigenvalue with multiplicity two, i.e. $\lambda_2=\lambda_1$ in (\ref{jump}), and in this case we should modify the image space in both components. Furthermore, very often in applications we have the eigenvalue $\lambda_1=-1$, due to $w_1=-1/2$, which leads to an image space of the type

$$Y_1=r_+\Lambda_-^{-s}\ell^{(0)}\{  \breve{H}^{0}(\R_+)\times L^2(\R_+)\}$$
with $\breve{H}^{0}(\R_+)=r_+\Lambda_-^{-1/2}H_+^{-1/2}$ being a proper dense subspace of $ L^2(\R_+)$ \cite{AmSpeckTeixeira98}. 

Finally, note that after we get the image normalization of a particular WHO in (\ref{equivalentWHO}), by means of the equivalence relation (\ref{equivalentrelation}), we achieve the image normalization of the operator ${\mathcal P}$ in (\ref{operatorP}).

%%%%%%%%%%%%%%%%%%%%%%%%%

\section{Boundary-transmission problems of higher order}  \label{example1}

From this section on, we analyze several examples of boundary-transmission conditions,  less general  than  (\ref{boundary1})-(\ref{boundary2}), but  still very significant from the applications point of view. First of all we retain only the higher order terms in the boundary-transmission conditions (\ref{boundary1})-(\ref{boundary2}) and such that they do not contain derivatives of mixed type. This assumption is also mathematically consistent with the fact that these terms fully describe the behavior at infinity of the Fourier symbol. Let us then consider the Problem ${\mathcal P}$ with the following higher order boundary-transmission conditions: order $m=(m_1,m_2)$ on the left half-line  and order $m'=(m_3,m_4)$ on the right half-line \footnote{In this section we intentionally use the notation of orders $m_j$, $j=1,2,3,4$ corresponding to the indices of the coefficients in order not to overload the formulas. We also simplify the index notations of the coefficients.}

\begin{equation} \label{orderm1,2}
 \begin{array}{ccl}
\left\{
\begin{array}{lll}
a_1^+{\varphi}_{m_1}^{+}+a_1^-{\varphi}_{m_1}^{-} +\check{a}_1^+\check{\varphi}_{m_1}^{+}+\check{a}_1^-\check{\varphi}_{m_1}^{-}&=& h_1\\
a_2^+{\varphi}_{m_2}^{+}+a_2^-{\varphi}_{m_2}^{-} +\check{a}_2^+\check{\varphi}_{m_2}^{+}+\check{a}_2^-\check{\varphi}_{m_2}^{-}&=& h_2
\end{array}
\right. \;\,\quad &\mbox{on}& {\R_-},
\end{array}
\end{equation}

\begin{equation} \label{orderm3,4}
 \begin{array}{ccl}
\left\{
\begin{array}{lll}
b_3^+{\varphi}_{m_3}^{+}+b_3^-{\varphi}_{m_3}^{-} +\check{b}_3^+\check{\varphi}_{m_3}^{+}+\check{b}_3^-\check{\varphi}_{m_3}^{-} &=& g_1\\
b_4^+{\varphi}_{m_4}^{+}+b_4^-{\varphi}_{m_4}^{-} +\check{b}_4^+\check{\varphi}_{m_4}^{+}+\check{b}_4^-\check{\varphi}_{m_4}^{-}&=& g_2
\end{array}
\right. \;\,\quad &\mbox{on}& {\R_+}
\end{array}
\end{equation}
where, at least for now, all four orders $m_j$, $j=1,2,3,4$ are supposed to be different, the notation ${\varphi}_{m_j}^{\pm}$ and $\check{\varphi}_{m_j}^{\pm}$ stands for the traces of the normal and tangential derivatives of order $m_j$, respectively, with  ${a}_{j}^{\pm}$,  ${b}_{j}^{\pm}$,  $\check{a}_{j}^{\pm}$,  $\check{b}_{j}^{\pm}$ denoting the corresponding coefficients. From the trace theorem and the representation formula, we conclude that $h=(h_1, h_2) \in H^{1/2-m_1}(\R_-)\times H^{1/2-m_2}(\R_-)$ and $g=(g_1,g_2)\in H^{1/2-m_3}(\R_+)\times H^{1/2-m_4}(\R_+)$.

The following theorem holds for the operator ${\mathcal P}$ associated with this boundary-transmission problem and the equivalent WHO, and is a direct consequence of Theorems \ref{traceoperator} and \ref{equivalencetoWHO}.

\begin{teo} \label{equivalencetoWHO1}
Let ${B}_{-}= {\mathcal F}^{-1}
{\Phi_{-}} \cdot {\mathcal F}$ and  $B_{+}= {\mathcal F}^{-1}
{\Phi_{+}} \cdot {\mathcal F}$  be the following linear bounded operators
$$
{B}_{-}: [H^{1/2}]^2 \rightarrow  H^{1/2-m_1}\times H^{1/2-m_2}
$$
$$
B_{+}: [H^{1/2}]^2 \rightarrow  H^{1/2-m_3}\times H^{1/2-m_4}
$$
with the non-degenerated Fourier symbols
\begin{equation} \label{Fourier12}
\Phi_-=
\left[
   \begin{array}{cc}
(-1)^{m_1}a_1^+{\beta}^{m_1}+\check{a}_1^+(-i\xi)^{m_1}
  & a_1^-{\beta}^{m_1}+\check{a}_1^-(-i\xi)^{m_1}  \\  
(-1)^{m_2}a_2^+{\beta}^{m_2}+\check{a}_2^+(-i\xi)^{m_2} & 
a_2^-{\beta}^{m_2}+\check{a}_2^-(-i\xi)^{m_2}  \end{array}
 \right]  
\end{equation}
and 
\begin{equation} \label{Fourier34}
\Phi_+=
\left[
   \begin{array}{cc}
(-1)^{m_3}b_3^+{\beta}^{m_3}+\check{b}_3^+(-i\xi)^{m_3}
  & b_3^-{\beta}^{m_3}+\check{b}_3^-(-i\xi)^{m_3}  \\  
(-1)^{m_4}b_4^+{\beta}^{m_4}+\check{b}_4^+(-i\xi)^{m_4} & 
b_4^-{\beta}^{m_4}+\check{b}_4^-(-i\xi)^{m_4}  \end{array}
 \right]  ,
\end{equation}
respectively, i.e.  $\rm {det}\Phi_+\Phi_-^{-1}\neq 0$. Moreover consider the restricted operator $\tilde{B}_-={\rm {Rst}} B_-: [H^{1/2}]^2 \rightarrow  H_+^{1/2-m_1}\times H_+^{1/2-m_2}$ with the Fourier symbol also given by (\ref{Fourier12}).
Then, the operator ${\mathcal P}$ given by
$$
\begin{array}{ccl}
{\mathcal P}: D({\mathcal P})& \rightarrow& H^{1/2-m_3} (\R_+)\times H^{1/2-m_4} (\R_+)\\
 \varphi & \rightarrow& {\mathcal P}\varphi=g,
\end{array}
$$
with $D({\mathcal P})$ defined as a closed subspace of  $H^1(\R\times \R_+)\bigoplus H^1(\R\times \R_-)$ and the solution $\varphi$ in (\ref{represent}) with traces $\varphi_0=({\varphi}_0^+,{\varphi}_0^- )= \tilde{B}_-^{-1}\left(v^{+}+\ell ^{(c)} h\right)$,  $v^{+}\in  H_+^{1/2-m_1}\times H_+^{1/2-m_2}$, is equivalent to the WHO
$$
\begin{array}{ccl}
W  : H^{1/2-m_1}_+ \times H^{1/2-m_2}_+& \rightarrow& H^{1/2-m_3}(\R_+) \times H^{1/2-m_4}(\R_+)\\
v^{+}& \rightarrow& Wv^{+}=g-r_{+}B_+\tilde{B}_-^{-1} \ell ^{(c)} h,
\end{array}
$$
i.e. $W=r_{+} {\mathcal F}^{-1}
{\Phi}\cdot {\mathcal F}$ with Fourier symbol $\Phi=\Phi_+\Phi_-^{-1}$. The equivalence relation is given by ${\mathcal P}=W\tilde{B}_- T_0$.
\end{teo}

A straightforward computation leads to the Fourier symbol of the equivalent WHO of the form

\begin{equation} \label{Fourierhigher}
\Phi=\frac{1}{{\rm det}\Phi_-}
\left[
   \begin{array}{cc}
A_{11}
  &A_{12} \\  
A_{21} & 
A_{22} \end{array}
 \right] ,
\end{equation}
with $${\rm det}\Phi_-=\beta^{m_1+m_2}((-1)^{m_1}a_1^+a_2^--(-1)^{m_2}a_2^+a_1^-)+
(-i\xi)^{m_1+m_2}(\check{a}_1^+\check{a}_2^--\check{a}_1^-\check{a}_2^+)+$$
$$(-i\xi)^{m_1}\beta^{m_2}((-1)^{m_1}\check{a}_1^+a_2^--
(-1)^{m_2}\check{a}_1^-a_2^+)+
(-i\xi)^{m_2}\beta^{m_1}((-1)^{m_1}a_1^+\check{a}_2^--(-1)^{m_2}a_1^-\check{a}_2^+),
$$ and entries
$$A_{11}=\beta^{m_2+m_3}((-1)^{m_3}a_2^-b_3^+-(-1)^{m_2}a_2^+b_3^-+(\check{a}_2^-\check{b}_3^++\check{a}_2^+\check{b}_3^-)(-i\xi\beta^{-1})^{m_2+m_3}+
$$
$$
((-1)^{m_3}\check{a}_2^-{b}_3^++\check{a}_2^+{b}_3^-)(-i\xi\beta^{-1})^{m_2}+
({a}_2^-\check{b}_3^+-(-1)^{m_2}{a}_2^+\check{b}_3^-)(-i\xi\beta^{-1})^{m_3}),
$$
$$A_{12}=\beta^{m_1+m_3}((-1)^{m_1}a_1^+b_3^--(-1)^{m_3}a_1^-b_3^++(\check{a}_1^+\check{b}_3^--{a}_1^-\check{b}_3^+)(-i\xi\beta^{-1})^{m_1+m_3}+
$$
$$
((-1)^{m_3}\check{a}_1^-{b}_3^++\check{a}_1^+{b}_3^-)(-i\xi\beta^{-1})^{m_1}+
((-1)^{m_1}{a}_1^+\check{b}_3^--\check{a}_1^-\check{b}_3^+)(-i\xi\beta^{-1})^{m_3}),
$$
$$A_{21}=\beta^{m_2+m_4}((-1)^{m_4}a_2^-b_4^+-(-1)^{m_2}a_2^+b_4^-+(\check{a}_2^-\check{b}_4^++\check{a}_2^+\check{b}_4^-)(-i\xi\beta^{-1})^{m_2+m_4}+
$$
$$
((-1)^{m_4}\check{a}_2^-{b}_4^++\check{a}_2^+{b}_4^-)(-i\xi\beta^{-1})^{m_2}+
({a}_2^-\check{b}_4^+-(-1)^{m_2}{a}_2^+\check{b}_4^-)(-i\xi\beta^{-1})^{m_4}),
$$
$$A_{22}=\beta^{m_1+m_4}((-1)^{m_1}a_1^+b_4^--(-1)^{m_4}a_1^-b_4^++(\check{a}_1^+\check{b}_4^--\check{a}_1^-\check{b}_4^+)(-i\xi\beta^{-1})^{m_1+m_4}+
$$
$$
(\check{a}_1^+{b}_4^--(-1)^{m_4}\check{a}_1^-{b}_4^+)(-i\xi\beta^{-1})^{m_1}+((-1)^{m_1}{a}_1^+\check{b}_4^--{a}_1^-\check{b}_4^+)(-i\xi\beta^{-1})^{m_4}),
$$
where $\beta(\xi)=\sqrt{\xi^2-k_0^2}$. The corresponding lifted Fourier symbol can be obtained based on the standard lifting procedure, i.e. taking $\Phi_{0} ={\rm diag}(\lambda_-^{1/2-m_3},\lambda_-^{1/2-m_4})\Phi {\rm diag}(\lambda_+^{m_1-1/2},\lambda_+^{m_2-1/2})$. We  obtain explicitely

\begin{equation} \label{Fourierhigherlifted}
\Phi_0=\frac{\rho}{{\rm det}\Phi_-}
\left[
   \begin{array}{cc}
\frac{(\xi+k_0)^{m_1}}{(\xi-k_0)^{m_3}}A_{11}
  &\frac{(\xi+k_0)^{m_2}}{(\xi-k_0)^{m_3}}A_{12} \\  
\frac{(\xi+k_0)^{m_1}}{(\xi-k_0)^{m_4}}A_{21} & 
\frac{(\xi+k_0)^{m_2}}{(\xi-k_0)^{m_4}}A_{22} \end{array}
 \right],
\end{equation}
where $\rho(\xi)=\sqrt{\frac{\xi-k_0}{\xi+k_0}}$.

In general, the given operator ${\mathcal P}$ is not normally solvable for arbitrary orders $m_j$ and coefficients ${a}_{j}^{\pm}$,  ${b}_{j}^{\pm}$,  $\check{a}_{j}^{\pm}$,  $\check{b}_{j}^{\pm}$, $j=1,2,3,4$.  For instance, if  all orders are even, then the not normally solvability does not depend on the coefficients. Although it will be very cumbersome to enumerate all the possibilities, the following easy to prove case shows the efficacy of the image normalization technique.  

\begin{teo} \label{allevencase}
Let $m_1=m_2=m_3=m_4=m\in 2\N_0$, i.e. all orders are equal to $m$, and $m$ is zero or an even number.  Furthermore, let $\rm {det}\Phi_0\neq 0$ in (\ref{Fourierhigherlifted}). Then the associated operator ${\mathcal P}$, and consequently the equivalent WHO $W$, are not normally solvable.
In this case we consider the corresponding {\it image normalized} operator $\breve{W}$ defined by

$$\breve{W}={\rm Rst}W: [H_+^{1/2-m}]^2\rightarrow  Y_1=r_+\Lambda_-^{-s}\ell^{(0)}\{\breve{H}^{0}(\R_+)\times \breve{H}^{0}(\R_+)\},$$
and $s=(1/2-m,1/2-m)$.
The image space $Y_1$ of the restricted operator $\breve{W}$ solves the {\it normalization problem} for the WHO, consequently  we look for solutions of ${\mathcal P}\varphi=g$, ${\mathcal P}: D({\mathcal P}) \rightarrow Y_1$, for which $g-r_+B_+\tilde{B}_-^{-1}\ell^{(c)}h \in Y_1$.
\end{teo}

\begin{pro}
The lifted Fourier symbol in (\ref{Fourierhigherlifted}) for $m=m_1=m_2=m_3=m_4\in 2\N_0$ simplifies to
$$\Phi_{0} =
\frac{\rho}{A}
\left[
   \begin{array}{cc}
B_{11}
  &B_{12} \\  
B_{21} & 
B_{22} \end{array}
 \right] 
$$
where
$$A=a_1^+a_2^--a_2^+a_1^-+
(\xi\beta^{-1})^{2m}(\check{a}_1^+\check{a}_2^--\check{a}_1^-\check{a}_2^+)+
(-i)^m(\xi\beta^{-1})^{m}(\check{a}_1^+a_2^--
\check{a}_1^-a_2^++a_1^+\check{a}_2^--a_1^-\check{a}_2^+),
$$ 
$$B_{11}=a_2^-b_3^+-a_2^+b_3^-+
(\xi\beta^{-1})^{2m}
(\check{a}_2^-\check{b}_3^++\check{a}_2^+\check{b}_3^-)+
(-i)^m(\xi\beta^{-1})^{m}
(\check{a}_2^-{b}_3^++\check{a}_2^+{b}_3^-+{a}_2^-\check{b}_3^+-{a}_2^+\check{b}_3^-),
$$
$$B_{12}=a_1^+b_3^--a_1^-b_3^++(\xi\beta^{-1})^{2m}(\check{a}_1^+\check{b}_3^--{a}_1^-\check{b}_3^+)+(-i)^m(\xi\beta^{-1})^{m}(\check{a}_1^-{b}_3^++\check{a}_1^+{b}_3^-+{a}_1^+\check{b}_3^--\check{a}_1^-\check{b}_3^+),
$$
$$B_{21}=a_2^-b_4^+-a_2^+b_4^-+(\xi\beta^{-1})^{2m}(\check{a}_2^-\check{b}_4^++\check{a}_2^+\check{b}_4^-)+
(-i)^m(\xi\beta^{-1})^{m}(\check{a}_2^-{b}_4^++\check{a}_2^+{b}_4^-+{a}_2^-\check{b}_4^+-{a}_2^+\check{b}_4^-),
$$
$$B_{22}=a_1^+b_4^--a_1^-b_4^++(\xi\beta^{-1})^{2m}(\check{a}_1^+\check{b}_4^--\check{a}_1^-\check{b}_4^+)+
(-i)^m(\xi\beta^{-1})^{m}(\check{a}_1^+{b}_4^--\check{a}_1^-{b}_4^++{a}_1^+\check{b}_4^--{a}_1^-\check{b}_4^+).
$$

Recall that $\rho(\xi)$, as well as $\xi\beta(\xi)^{-1}$, tends to $\pm 1$ as $\xi$ tends to $\pm \infty$, respectively. But here the $\xi\beta(\xi)^{-1}$ factors are all raised to an even power: $2m$ or $m$. Thus $\Phi_0(-\infty)=-\Phi_0(+\infty)$ and for the Fredholm criterium (\ref{criterium}) one has
$$
\mu \Phi_0(-\infty)+(1-\mu) \Phi_0(+\infty)=(1-2\mu)\Phi_0(+\infty),
$$
which degenerates for $\mu=1/2$, i.e.  $\Phi_0$ doesn't fulfill the Fredholm criterium for $\mu=1/2$. After some calculations we arrive at
 $$
\Phi_0^{-1}(-\infty)\Phi_0(+\infty)=\left[
   \begin{array}{cc}
-1
  &0 \\  
0 & 
-1 \end{array}
 \right] 
 $$
Thus the result follows from Theorem \ref{Main Theorem-matrix}, since the jump at infinity (\ref{jump}) has a diagonal form with one eigenvalue  $\lambda=-1$ with multiplicity two.
\end{pro}

We remark once again that an analogous result can be obtained for even orders $m_j$ not necessarily  all equals, see e.g. Theorem \ref{normalizeallevencase}, only in that case the calculations are more complicated. 

%%%%%%%%%%%%%%%%%%%%%%%%%%%%%%%%%%%%%%%

\section{Boundary-transmission problems of pairwise normal type}  \label{example2}

We formulate now a particular case of boundary-transmission conditions of the form (\ref{orderm1,2})-(\ref{orderm3,4}), namely consider on both upper banks of $\R_-$ and $\R_+$ boundary-transmission conditions with normal derivatives of a given order, say $m_1$, and on both lower banks of $\R_-$ and $\R_+$ boundary-transmission conditions with normal derivatives of another order, say $m_2$. We can similarly to the previous Section \ref{example1}, define the associated operator ${\mathcal P}$ to the problem and study its normal solvability together with the equivalent WHO.  Let us consider, together with the Helmholtz equation (\ref{Helmholtz}),   the following boundary-transmission conditions of orders $m=(m_1,m_2)$ and $m'=(m_1,m_2)$.

\begin{equation} \label{normalm1m2}
 \begin{array}{ccl}
\left\{
\begin{array}{lll}
a_1^+{\varphi}_{m_1}^{+}+a_1^-{\varphi}_{m_1}^{-} &=& h_1\\
a_2^+{\varphi}_{m_2}^{+}+a_2^-{\varphi}_{m_2}^{-}&=& h_2
\end{array}
\right. \;\,\quad &\mbox{on}& {\R_-},
\end{array}
\end{equation}

\begin{equation} \label{normalm1andm2}
 \begin{array}{ccl}
\left\{
\begin{array}{lll}
b_1^+{\varphi}_{m_1}^{+}+b_1^-{\varphi}_{m_1}^{-}  &=& g_1\\
b_2^+{\varphi}_{m_2}^{+}+b_2^-{\varphi}_{m_2}^{-} &=& g_2
\end{array}
\right. \;\,\quad &\mbox{on}& {\R_+}
\end{array}
\end{equation}
where  $m_1\neq m_2$, $h=(h_1, h_2) \in H^{1/2-m_1}(\R_-)\times H^{1/2-m_2}(\R_-)$ and $g=(g_1,g_2) \in H^{1/2-m_1}(\R_+)\times H^{1/2-m_2}(\R_+)$.

Here we should consider two cases: when $m_1+m_2$ is even or zero, and when $m_1+m_2$ is odd, due to the following necessary and sufficient conditions for the operator ${\mathcal P}$ and the equivalent WHO be of normal type.

\begin{teo} \label{normaltypeWHO} 
Consider the associated operator ${\mathcal P}$ 
$$
\begin{array}{ccl}
{\mathcal P}: D({\mathcal P})& \rightarrow& H^{1/2-m_1} (\R_+)\times H^{1/2-m_2} (\R_+)\\
 \varphi & \rightarrow& {\mathcal P}\varphi=g,
\end{array}
$$
and the equivalent WHO
$$
W= r_{+} {\mathcal F}^{-1}
{\Phi}\cdot {\mathcal F} : H^{1/2-m_1}_+ \times H^{1/2-m_2}_+ \rightarrow H^{1/2-m_1}(\R_+) \times H^{1/2-m_2}(\R_+)
$$
with Fourier symbol
\begin{equation} \label{Fouriernormalm1m2}
\Phi=\frac{1}{A}
\left[
   \begin{array}{cc}
(-1)^{m_1}a_2^-b_1^+-(-1)^{m_2}a_2^+b_1^-
  &(-1)^{m_1}(a_1^+b_1^--a_1^-b_1^+)\beta^{m_1-m_2}  \\  
(-1)^{m_2}(a_2^-b_2^+-a_2^+b_2^-)\beta^{m_2-m_1}  & 
(-1)^{m_1}a_1^+b_2^--(-1)^{m_2}a_1^-b_2^+  \end{array}
 \right] , 
\end{equation}
where $A=(-1)^{m_1}a_1^+a_2^--(-1)^{m_2}a_1^-a_2^+$. Then the operator $W$, and consequently the operator ${\mathcal P}$, are of normal type iff
\begin{equation} \label{normalcondit}
\frac{a_1^+a_2^-b_1^+b_2^-+a_1^-a_2^+b_1^-b_2^+}{a_1^-a_2^+b_1^+b_2^-+a_1^+a_2^-b_1^-b_2^+}\neq (-1)^{m_1+m_2} 
\end{equation}
\end{teo}

\begin{pro}
The operator $W$ is obtained as in Section \ref{example1} from definition (\ref{equivalentWHO}), i.e. $W=r_{+}B_+\tilde{B}_-^{-1}$, where now the operators $\tilde{B}_-$ and ${B}_+$ have the folllowing Fourier symbols 
$$
\Phi_-=
\left[
   \begin{array}{cc}
a_1^+({-\beta})^{m_1}
  & a_1^-{\beta}^{m_1} \\  
a_2^+({-\beta})^{m_2} & 
a_2^-{\beta}^{m_2}  \end{array}
 \right]  
$$
$$
\Phi_+=
\left[
   \begin{array}{cc}
b_1^+({-\beta})^{m_1}
  & b_1^-{\beta}^{m_1} \\  
b_2^+({-\beta})^{m_2} & 
b_2^-{\beta}^{m_2}  \end{array}
 \right]. 
 $$
Note that these are particular cases of symbols (\ref{Fourier12}) and (\ref{Fourier34}) for zero coefficients of the tangential derivatives and $m_3=m_1$, $m_4=m_2$. Then from (\ref{Fouriernormalm1m2}) we can obtain the  lifted Fourier symbol doing $\Phi_{0} ={\rm diag}(\lambda_-^{1/2-m_1},\lambda_-^{1/2-m_2})\Phi {\rm diag}(\lambda_+^{m_1-1/2},\lambda_+^{m_2-1/2})$ 
\begin{equation} \label{Fouriernormalm1m2lifted}
\Phi_{0} =
\frac{1}{A}
\left[
   \begin{array}{cc}
\left((-1)^{m_1}a_2^-b_1^+-(-1)^{m_2}a_2^+b_1^- \right)\rho^{1-2m_1}
  &(-1)^{m_1}(a_1^+b_1^--a_1^-b_1^+) \rho^{1-m_1-m_2} \\  
(-1)^{m_2}(a_2^-b_2^+-a_2^+b_2^- ) \rho^{1-m_1-m_2}  & 
\left((-1)^{m_1}a_1^+b_2^--(-1)^{m_2}a_1^-b_2^+\right)\rho^{1-2m_2} \end{array}
 \right] , 
\end{equation}
where $A=(-1)^{m_1}a_1^+a_2^--(-1)^{m_2}a_1^-a_2^+$. Finally, condition ${\rm det}\Phi_0\neq 0$  is equivalent to
$$
   \begin{array}{cl}
a_1^+a_2^-b_1^+b_2^-+a_1^-a_2^+b_1^-b_2^++&\\
(-1)^{m_1+m_2+1}&\!\!\!\!\!\!\!\!\!\!\!\!\!\!\!\!\!\!\!\!\left( a_1^-a_2^-b_1^+b_2^++a_1^+a_2^+b_1^-b_2^-+(a_1^-b_1^+-a_1^+b_1^-)(a_2^-b_2^+-a_2^+b_2^-)  \right)\neq 0
\end{array}
$$
which can also be simplified to equation (\ref{normalcondit}).
\end{pro}

Condition (\ref{normalcondit}) means that for orders $m_1+m_2\in 2\N_0$  the operators ${\mathcal P}$ and $W$ are of normal type iff $a_1^+a_2^-\neq a_1^-a_2^+$ and $b_1^+b_2^-\neq b_1^-b_2^+$. On the other hand, if $m_1+m_2\in 2\N_0+1$, then the operators ${\mathcal P}$ and $W$ are of normal type iff $a_1^+a_2^-\neq -a_1^-a_2^+$ and $b_1^+b_2^-\neq -b_1^-b_2^+$. The first case gives place to the following theorem on the image normalization of $W$ and ${\mathcal P}$.

\begin{teo} \label{normalizevencase}
Let $m_1+m_2\in 2\N_0$ in the boundary-transmission conditions (\ref{normalm1m2})-(\ref{normalm1andm2}) and assume that (\ref{normalcondit}) holds. Then the operator ${\mathcal P}$ and the equivalent operator $W$ are not-normally solvable. In this case, the image space of the {\it image normalized} operator $\breve{W}$ given by

$$\breve{W}={\rm Rst}W: [H_+^{1/2-m_1}]^2\rightarrow  r_+\Lambda_-^{-s}\ell^{(0)}\{\breve{H}^{0}(\R_+)\times \breve{H}^{0}(\R_+)\},$$
with $s=(1/2-m_1,1/2-m_2)$,
solves the {\it normalization problem} for the WHO. The image normalization of operator ${\mathcal P}$ is achieved by substituting $W$ by $\breve{W}$ in the equivalence relation (\ref{equivalentrelation}).
\end{teo}

\begin{pro}
For  $m_1+m_2\in 2\N_0$ we have $(-1)^{m_1}=(-1)^{m_2}$, and the lifted Fourier symbol in (\ref{Fouriernormalm1m2lifted}) simplifies to
$$\Phi_{0} =
\frac{1}{a_1^+a_2^--a_1^-a_2^+}
\left[
   \begin{array}{cc}
\left(a_2^-b_1^+-a_2^+b_1^- \right)\rho^{1-2m_1}
  &\left(a_1^+b_1^--a_1^-b_1^+\right) \rho^{1-m_1-m_2} \\  
\left(a_2^-b_2^+-a_2^+b_2^- \right) \rho^{1-m_1-m_2}  & 
\left(a_1^+b_2^--a_1^-b_2^+\right)\rho^{1-2m_2} \end{array}
 \right].  
$$
Remark that $\rho(\xi)^{1-m_1-m_2}$, as well as $\rho(\xi)^{1-2m_j}$, $j=1,2$, tends to $\pm 1$ as $\xi$ tends to $\pm \infty$, respectively.  Thus $\Phi_0(-\infty)=-\Phi_0(+\infty)$ and the Fredholm criterium (\ref{criterium}) gives
$$
\mu \Phi_0(-\infty)+(1-\mu) \Phi_0(+\infty)=(1-2\mu)\Phi_0(+\infty),
$$
which degenerates for $\mu=1/2$, i.e.  $\Phi_0$ doesn't fulfill the Fredholm criterium for $\mu=1/2$. We now arrive at
 $$
\Phi_0^{-1}(-\infty)\Phi_0(+\infty)=
\left[
   \begin{array}{cc}
\frac{(a_1^+a_2^--a_1^-a_2^+)(b_1^+b_2^--b_1^-b_2^+)}{(a_1^-a_2^+-a_1^+a_2^-)(b_1^+b_2^--b_1^-b_2^+)}
  &0 \\  
0 & 
\frac{(a_1^+a_2^--a_1^-a_2^+)(b_1^+b_2^--b_1^-b_2^+)}{(a_1^-a_2^+-a_1^+a_2^-)(b_1^+b_2^--b_1^-b_2^+)}
 \end{array}
 \right] =
\left[
   \begin{array}{cc}
-1
  &0 \\  
0 & 
-1 \end{array}
 \right]. 
 $$
Therefore the result is a consequence of the Theorem \ref{Main Theorem-matrix}, applied to the jump at infinity of diagonal form with one eigenvalue  $\lambda=-1$ with multiplicity two.
\end{pro}

%%%%%%%%%%%%%%%%%

For orders $m_1+m_2\in 2\N_0+1$ the following normalization theorem shows us  that the image normalization can also depend on the coefficients.

\begin{teo} \label{normalizeoddcase}
Let $m_1+m_2\in 2\N_0+1$ in the boundary-transmission conditions (\ref{normalm1m2})-(\ref{normalm1andm2}) and assume that (\ref{normalcondit}) holds. Then the operator ${\mathcal P}$ and the equivalent operator $W$ are not-normally solvable iff there exists a solution $\theta \in [-1,1]$ for the equation
\begin{equation} \label{normalsolvablecondit}
\theta^2=\frac{(a_1^+b_1^--a_1^-b_1^+)(a_2^+b_2^--a_2^-b_2^+)}{(a_1^+b_2^-+a_1^-b_2^+)(a_2^-b_1^++a_2^+b_1-)} .
\end{equation}
In this case,  the image space of the {\it image normalized} operator $\breve{W}$ given by

$$\breve{W}={\rm Rst}W: [H_+^{1/2-m_1}]^2\rightarrow  r_+\Lambda_-^{-s}T\ell^{(0)}\{\breve{H}^{-i\tau}(\R_+)\times L_2(\R_+)\},$$
with $s=(1/2-m_1,1/2-m_2)$ solves the {\it normalization problem} for the WHO. Here $T$ is the matrix  which allows  the diagonalization $\Phi_0^{-1}(-\infty)\Phi_0(+\infty)=T^{-1}{\rm diag}(\lambda_1,\lambda_2)T$ for which the eigenvalue $\lambda_1$ has argument equal to $-\pi$ and $\tau=\frac{1}{2\pi}\log \left| \frac{\Phi_0(-\infty)}{\Phi_0(+\infty)} \right|$ corresponds to $\lambda_1$. Moreover, the image normalization of operator ${\mathcal P}$ is achieved by substituting $W$ by $\breve{W}$ in the equivalence relation (\ref{equivalentrelation}).
\end{teo}

\begin{pro}
For  $m_1+m_2\in 2\N_0+1$ we have $(-1)^{m_1}=-(-1)^{m_2}$, and the lifted Fourier symbol in (\ref{Fouriernormalm1m2lifted}) simplifies to
$$\Phi_{0} =
\frac{1}{a_1^+a_2^-+a_1^-a_2^+}
\left[
   \begin{array}{cc}
\left(a_2^-b_1^++a_2^+b_1^- \right)\rho^{1-2m_1}
  &\left(a_1^+b_1^--a_1^-b_1^+\right) \rho^{1-m_1-m_2} \\  
\left(a_2^+b_2^--a_2^-b_2^+\right) \rho^{1-m_1-m_2}  & 
\left(a_1^+b_2^-+a_1^-b_2^+\right)\rho^{1-2m_2} \end{array}
 \right].  
$$
Remark  that now, while $\rho(\xi)^{1-m_1-m_2}$ tends to one, $\rho(\xi)^{1-2m_j}$, $j=1,2$, tends to $\pm 1$ as $\xi$ tends to $\pm \infty$, respectively.  Therefore for the Fredholm criterium (\ref{criterium}) we have
$$
\mu \Phi_0(-\infty)+(1-\mu) \Phi_0(+\infty)=\frac{1}{a_1^+a_2^-+a_1^-a_2^+}
\left[
   \begin{array}{cc}
(1-2\mu)\left(a_2^-b_1^++a_2^+b_1^- \right)
  &a_1^+b_1^--a_1^-b_1^+ \\  
a_2^+b_2^--a_2^-b_2^+ & 
(1-2\mu)\left(a_1^+b_2^-+a_1^-b_2^+\right) \end{array}
 \right],
$$
which degenerates for
$$(1-2\mu)^2(a_2^-b_1^++a_2^+b_1^-)(a_1^+b_2^-+a_1^-b_2^+)-(a_1^+b_1^--a_1^-b_1^+ )(a_2^+b_2^--a_2^-b_2^+)=0$$
or equivalently when (\ref{normalsolvablecondit}) holds, where we introduced $\theta=1-2\mu$.

Since one has
$$ \Phi_0(-\infty)=\frac{1}{b_1^+b_2^-+b_1^-b_2^+}
\left[
   \begin{array}{cc}
-a_1^+b_2^--a_1^-b_2^+ 
  &-a_1^+b_1^-+a_1^-b_1^+ \\  
-a_2^+b_2^-+a_2^-b_2^+ & 
-a_2^+b_1^--a_2^-b_1^+ \end{array}
 \right],
$$
we get
 $$
\Phi_0^{-1}(-\infty)\Phi_0(+\infty)=
\left[
   \begin{array}{cc}
\frac{(\theta^2-1)(a_1^+b_2^-+a_1^-b_2^+)(a_2^-b_1^++a_2^+b_1^-)}{(a_1^+a_2^-+a_1^-a_2^+)(b_1^+b_2^-+b_1^-b_2^+)}
  &\frac{2(a_1^+b_2^-+a_1^-b_2^+)(a_1^-b_1^+-a_1^+b_1^-)}{(a_1^+a_2^-+a_1^-a_2^+)(b_1^+b_2^-+b_1^-b_2^+)} \\  
\frac{2(a_2^-b_1^++a_2^+b_1^-)(a_2^-b_2^+-a_2^+b_2^-)}{(a_1^+a_2^-+a_1^-a_2^+)(b_1^+b_2^-+b_1^-b_2^+)} & 
\frac{(1-\theta^2)(a_1^+b_2^-+a_1^-b_2^+)(a_2^-b_1^++a_2^+b_1^-)}{(a_1^+a_2^-+a_1^-a_2^+)(b_1^+b_2^-+b_1^-b_2^+)}
 \end{array}
 \right] 
 $$
 or, after introducing the notations
 $$A_{11}=\frac{(a_1^+b_2^-+a_1^-b_2^+)(a_2^-b_1^++a_2^+b_1^-)}{(a_1^+a_2^-+a_1^-a_2^+)(b_1^+b_2^-+b_1^-b_2^+)},
$$ 
$$A_{12}=\frac{2(a_1^+b_2^-+a_1^-b_2^+)(a_1^-b_1^+-a_1^+b_1^-)}{(a_1^+a_2^-+a_1^-a_2^+)(b_1^+b_2^-+b_1^-b_2^+)},
$$
$$A_{21}=\frac{2(a_2^-b_1^++a_2^+b_1^-)(a_2^-b_2^+-a_2^+b_2^-)}{(a_1^+a_2^-+a_1^-a_2^+)(b_1^+b_2^-+b_1^-b_2^+)},$$
we arrive at
 $$
\Phi_0^{-1}(-\infty)\Phi_0(+\infty)=
\left[
   \begin{array}{cc}
(\theta^2-1)A_{11}
  &A_{12} \\  
A_{21}& 
(1-\theta^2)A_{11}
 \end{array}
 \right]. 
 $$
The eigenvalues of the jump at infinity matrix are
\begin{equation} \label{eigenvalues}
 \lambda_j=\pm \sqrt {(1-\theta^2)^2A_{11}^2+A_{12}A_{21}},\,\,j=1,2 ,
\end{equation}
and the result is a consequence of the Theorem \ref{Main Theorem-matrix},  with the choice of $T$ to be the matrix, possible with a permutation of columns, that allows the eigenvalue $\lambda_1$ in (\ref{eigenvalues}) to be the one with argument equal to  $-\pi$.
\end{pro}

%%%%%%%%%%%%%%%%%%%%%%%%%%%%%%%%%

\section{Boundary-transmission problems with oblique derivatives}  \label{example3}

Finally, we analyze a particular case of boundary-transmission conditions of the form (\ref{orderm1,2})-(\ref{orderm3,4}), when we have a boundary condition with normal derivatives of order $m_1$ on both banks of $\R_-$ , and boundary-transmission conditions with normal and tangential derivatives of  order $m_2$ on both banks of $\R_+$.  Let us consider, together with the Helmholtz equation (\ref{Helmholtz}), the following boundary-transmission conditions of orders $m=(m_1,m_1)$ on the left half-line  and order $m'=(m_2,m_2)$ on the right half-line.

\begin{equation} \label{normalm1}
 \begin{array}{ccl}
\left\{
\begin{array}{lll}
{\varphi}_{m_1}^{-} &=& h_1\\
{\varphi}_{m_1}^{-}&=& h_2
\end{array}
\right. \;\,\quad &\mbox{on}& {\R_-},
\end{array}
\end{equation}

\begin{equation} \label{obliquem2}
 \begin{array}{ccl}
\left\{
\begin{array}{lll}
b_2^+{\varphi}_{m_2}^{+}+b_2^-{\varphi}_{m_2}^{-} +\check{b}_2^+\check{\varphi}_{m_2}^{+}+\check{b}_2^-\check{\varphi}_{m_2}^{-} &=& g_1\\
c_2^+{\varphi}_{m_2}^{+}+c_2^-{\varphi}_{m_2}^{-} +\check{c}_2^+\check{\varphi}_{m_2}^{+}+\check{c}_2^-\check{\varphi}_{m_2}^{-} &=& g_2\\
\end{array}
\right. \;\,\quad &\mbox{on}& {\R_+}
\end{array}
\end{equation}
where in general  $m_1\neq m_2$, $h=(h_1, h_2) \in [H^{1/2-m_1}(\R_-)]^2$ and $g=(g_1,g_2) \in [H^{1/2-m_2}(\R_+)]^2$.

In this section we should consider different cases depending on the parity of the orders $m_j$, $j=1,2$. The following necessary and sufficient conditions hold for the operator ${\mathcal P}$ and the equivalent WHO be of normal type.

\begin{teo} \label{normaltypeWHO} 
Consider the associated operator ${\mathcal P}$ 
$$
\begin{array}{ccl}
{\mathcal P}: D({\mathcal P})& \rightarrow& [H^{1/2-m_2}(\R_+)]^2\\
 \varphi & \rightarrow& {\mathcal P}\varphi=g,
\end{array}
$$
and the equivalent WHO
$$
W = r_{+} {\mathcal F}^{-1}
{\Phi}\cdot {\mathcal F}: [H_+^{1/2-m_1}]^2 \rightarrow [H^{1/2-m_2}(\R_+)]^2
$$
with Fourier symbol
\begin{equation} \label{Fourierobliquem1m2}
\Phi=(-1)^{m_1}\beta^{m_1-m_2}
\left[
   \begin{array}{cc}
(-1)^{m_2}b_2^++\check{b}_2^+(-i\xi\beta^{-1})^{m_2}
  &(-1)^{m_1}\left(b_2^-+\check{b}_2^-(-i\xi\beta^{-1})^{m_2}\right)  \\  
(-1)^{m_2}c_2^++\check{c}_2^+(-i\xi\beta^{-1})^{m_2} & 
(-1)^{m_1}\left(c_2^-+\check{c}_2^-(-i\xi\beta^{-1})^{m_2}\right)
  \end{array}
 \right]  
\end{equation}
Then both operators $W$ and  ${\mathcal P}$ are of normal type iff
\begin{equation} \label{normalcondition}
   \begin{array}{cl}
(-1)^{m_2}(b_2^+c_2^--b_2^-c_2^+)+&\!\!\left((-1)^{m_2}(b_2^+\check{c}_2^--\check{b}_2^-c_2^+)+\check{b}_2^+c_2^--b_2^-\check{c}_2^+)\right)(-i\xi\beta^{-1})^{m_2}\\
\,\,\,&+(\check{b}_2^+\check{c}_2^--\check{b}_2^-\check{c}_2^+)(-i\xi\beta^{-1})^{2m_2}\neq 0.
   \end{array}
\end{equation}
\end{teo}

\begin{pro}
As before, we obtain first the Fourier symbol in (\ref{Fourierobliquem1m2}) from the Fourier symbols of the operators $\tilde{B}_-$ and $B_+$, and then come to the lifted symbol 
\begin{equation} \label{Fourierobliqueliftedm1m2}
\Phi_0=(-1)^{m_1}\rho^{1-m_1-m_2}
\left[
   \begin{array}{cc}
(-1)^{m_2}b_2^++\check{b}_2^+(-i\xi\beta^{-1})^{m_2}
  &(-1)^{m_1}\left(b_2^-+\check{b}_2^-(-i\xi\beta^{-1})^{m_2}\right)  \\  
(-1)^{m_2}c_2^++\check{c}_2^+(-i\xi\beta^{-1})^{m_2} & 
(-1)^{m_1}\left(c_2^-+\check{c}_2^-(-i\xi\beta^{-1})^{m_2}\right)
  \end{array}
 \right] . 
\end{equation}
 Thus, condition (\ref{normalcondition}) follows from the assumption that ${\rm det}\Phi_0\neq 0$.  
\end{pro}

We must consider now four situations: both orders $m_j$ are zero or even, both orders are odd, $m_1$ is zero or even and $m_2$ is odd, and the way around, $m_1$ is odd and $m_2$ is zero or even. These four cases give rise to the following four results.

\begin{teo} \label{normalizeallevencase}
Let $m_1, m_2\in 2\N_0$ in the boundary-transmission conditions (\ref{normalm1})-(\ref{obliquem2}) and assume that (\ref{normalcondition}) holds. Then the operator ${\mathcal P}$ and the equivalent operator $W$ are not-normally solvable. In this case, the image space of the {\it image normalized} operator $\breve{W}$ defined by

$$\breve{W}={\rm Rst}W: [H_+^{1/2-m_1}]^2\rightarrow  r_+\Lambda_-^{-s}\ell^{(0)}\{\breve{H}^{0}(\R_+)\times \breve{H}^{0}(\R_+)\},$$
with $s=(1/2-m_2,1/2-m_2)$,
solves the {\it normalization problem} for the WHO. Furthermore, the image normalization of operator ${\mathcal P}$ is achieved by substituting $W$ by $\breve{W}$ in the equivalence relation (\ref{equivalentrelation}).
\end{teo}

\begin{pro}
For $m_1=m_2$ this is a direct consequence of Theorem \ref{allevencase}. For different orders $m_1\neq m_2$ we also arrive at $\Phi_0(-\infty)=-\Phi_0(+\infty)$, since the lifted Fourier symbols simplifies to 
$$\Phi_{0} =\rho^{1-m_1-m_2}
\left[
   \begin{array}{cc}
b_2^++i^{m_2}\check{b}_2^+(\xi\beta^{-1})^{m_2}
  &b_2^-+i^{m_2}\check{b}_2^-(\xi\beta^{-1})^{m_2}  \\  
c_2^++i^{m_2}\check{c}_2^+(\xi\beta^{-1})^{m_2} & 
c_2^-+i^{m_2}\check{c}_2^-(\xi\beta^{-1})^{m_2}
\end{array}
 \right].  
$$
It is not difficult to see that we obtain here for the jump at infinity, once more, a diagonal matrix with $-1$ in the diagonal entries.
\end{pro}

%%%%%%%%%%%%%%%%%

\begin{teo} \label{normalizealloddcase}
Let $m_1, m_2\in 2\N_0+1$ in the boundary-transmission conditions (\ref{normalm1})-(\ref{obliquem2}) and assume that (\ref{normalcondition}) holds. Then the operators ${\mathcal P}$ and $W$ are not-normally solvable if there exists a solution $\theta \in [-1,1]$ for the equation
\begin{equation} \label{normalsolvablecondition}
(b_2^+c_2^--b_2^-c_2^+)\theta^2+i^{m_2}(\check{b}_2^+c_2^--b_2^+\check{c}_2^-+\check{b}_2^-c_2^+-b_2^-\check{c}_2^+)\theta+\check{b}_2^+\check{c}_2^--\check{b}_2^-\check{c}_2^+  =0 .
\end{equation}
In this case,  the image space of the {\it image normalized} operator $\breve{W}$ given by

$$\breve{W}={\rm Rst}W: [H_+^{1/2-m_1}]^2\rightarrow  r_+\Lambda_-^{-s}T\ell^{(0)}\{\breve{H}^{-i\tau}(\R_+)\times L_2(\R_+)\},$$
with $s=(1/2-m_2,1/2-m_2)$ solves the {\it normalization problem} for the WHO. The matrix $T$ is chosen to be the matrix  in the diagonalization $\Phi_0^{-1}(-\infty)\Phi_0(+\infty)=T^{-1}{\rm diag}(\lambda_1,\lambda_2)T$ for which the eigenvalue $\lambda_1$ has argument equal to $-\pi$ and  $\tau=\frac{1}{2\pi}\log \left| \frac{\Phi_0(-\infty)}{\Phi_0(+\infty)} \right|$ corresponds to the eigenvalue $\lambda_1$. The image normalization of the operator ${\mathcal P}$ is achieved by substituting $W$ by $\breve{W}$ in the equivalence relation (\ref{equivalentrelation}).
\end{teo}

\begin{pro}
For  $m_1+m_2\in 2\N_0+1$ we have $(-1)^{m_1}=(-1)^{m_2}=-1$, and the lifted Fourier symbol in (\ref{Fourierobliqueliftedm1m2}) now simplifies to
$$\Phi_{0} =\rho^{1-m_1-m_2}
\left[
   \begin{array}{cc}
b_2^++i^{m_2}\check{b}_2^+(\xi\beta^{-1})^{m_2}
  &b_2^--i^{m_2}\check{b}_2^-(\xi\beta^{-1})^{m_2}  \\  
c_2^++i^{m_2}\check{c}_2^+(\xi\beta^{-1})^{m_2} & 
c_2^--i^{m_2}\check{c}_2^-(\xi\beta^{-1})^{m_2}
\end{array}
 \right].  
$$
Here, while $\rho(\xi)^{1-m_1-m_2}$ tends to one, $(\xi\beta(\xi)^{-1})^{m_2}$ tends to $\pm 1$ as $\xi$ tends to $\pm \infty$, respectively.  Thus the Fredholm criterium (\ref{criterium}) applied to the lifted Fourier symbol gives
$$\mu \Phi_0(-\infty)+(1-\mu) \Phi_0(+\infty)=
\left[
   \begin{array}{cc}
b_2^+\theta+i^{m_2}\check{b}_2^+
  &b_2^-\theta-i^{m_2}\check{b}_2^-  \\  
c_2^+\theta+i^{m_2}\check{c}_2^+& 
c_2^-\theta-i^{m_2}\check{c}_2^-
\end{array}
 \right] ,
$$
where we introduced the former notation $\theta=1-2\mu$. Remark that once  $i^{2m_2}=-1$,  the determinant of the last matrix equals zero when (\ref{normalsolvablecondition}) holds.

Since we have
 $$
\Phi_0^{-1}(-\infty)\Phi_0(+\infty)=
\left[
   \begin{array}{cc}
-\frac{A+i^{m_2}B}{C}
  &\frac{2i^{m_2}D}{C} \\  
\frac{2i^{m_2}E}{C} & 
\frac{-A+i^{m_2}B}{C}
 \end{array}
 \right],  
 $$
 where we introduced the notations
 $$A=b_2^-c_2^+-b_2^+c_2^-+\check{b}_2^+\check{c}_2^--\check{b}_2^-\check{c}_2^+,\,\,B=\check{b}_2^-\check{c}_2^++{b}_2^-\check{c}_2^+-\check{b}_2^+{c}_2^--{b}_2^+\check{c}_2^-,
$$ 
$$C=b_2^-c_2^+-b_2^+c_2^-+\check{b}_2^-\check{c}_2^-+\check{b}_2^-\check{c}_2^++i^{m_2} B,\,\,D={b}_2^-\check{c}_2^--\check{b}_2^-{c}_2^-,\,\,E={b}_2^+\check{c}_2^+-\check{b}_2^+{c}_2^+,
$$
we get the eigenvalues
\begin{equation} \label{eigenvaluesodd}
 \lambda_j=\frac{-A\pm i \sqrt {B^2+4DE}}{C},\,\,j=1,2.
\end{equation}
Therefore the result is a consequence of the Theorem \ref{Main Theorem-matrix}, applied to the jump at infinity with the choice of the eigenvalue $\lambda_1$ in (\ref{eigenvaluesodd}) to be the one for which the argument is equal to $-\pi$.
\end{pro}

%%%%%%%%%%%%%%%%%%%%%%%%%

\begin{teo} \label{normalizeevenoddcase}
Let $m_1\in 2\N_0$ and $m_2\in 2\N_0+1$ in the boundary-transmission conditions (\ref{normalm1})-(\ref{obliquem2}) and assume that (\ref{normalcondition}) holds. Then the operator ${\mathcal P}$ and the equivalent operator $W$ are not-normally solvable if there exists a solution $\theta \in [-1,1]$ for the equation
\begin{equation} \label{normalsolvablecondition1}
(\check{b}_2^-\check{c}_2^+-\check{b}_2^+\check{c}_2^-)\theta^2-i^{m_2}(\check{b}_2^+c_2^--\check{c}_2^-b_2^++\check{b}_2^-c_2^+-b_2^-\check{c}_2^+)\theta+b_2^-c_2^+- b_2^+c_2^- =0 .
\end{equation}
In this case,  the image space of the {\it image normalized} operator $\breve{W}$ given by

$$\breve{W}={\rm Rst}W: [H_+^{1/2-m_1}]^2\rightarrow  r_+\Lambda_-^{-s}T\ell^{(0)}\{\breve{H}^{-i\tau}(\R_+)\times L_2(\R_+)\},$$
with $s=(1/2-m_2,1/2-m_2)$ solves the {\it normalization problem} for the WHO. The matrix $T$ is chosen to be the matrix  in the diagonalization $\Phi_0^{-1}(-\infty)\Phi_0(+\infty)=T^{-1}{\rm diag}(\lambda_1,\lambda_2)T$ for which the eigenvalue $\lambda_1$ has argument equal to $-\pi$ and  $\tau=\frac{1}{2\pi}\log \left| \frac{\Phi_0(-\infty)}{\Phi_0(+\infty)} \right|$ corresponds to the eigenvalue $\lambda_1$. The image normalization of the operator ${\mathcal P}$ can be achieved by substituting $W$ by $\breve{W}$ in the equivalence relation (\ref{equivalentrelation}).
\end{teo}

\begin{pro}
The proof follows the same steps as the proof of Theorem \ref{normalizealloddcase}. In this case the jump at infinity is characterized by the matrix
 $$
\Phi_0^{-1}(-\infty)\Phi_0(+\infty)=
\left[
   \begin{array}{cc}
\frac{A+i^{m_2}B}{C}
  &\frac{2i^{m_2}D}{C} \\  
\frac{2i^{m_2}E}{C} & 
\frac{A-i^{m_2}B}{C}
 \end{array}
 \right],  
 $$
where the notations are the same as the ones used in the proof of Theorem \ref{normalizealloddcase}. The eigenvalues are given by
  $$
 \lambda_j=\frac{A\pm i \sqrt {B^2+4DE}}{C},\,\,\,j=1,2
 $$
and we must choose the diagonalization which gives the value of $-\pi$ for the argument of $\lambda_1$.
\end{pro}

%%%%%%%%%%%%%%%%%%%%%%%%%

\begin{teo} \label{normalizeoddevencase}
Let $m_1\in 2\N_0+1$ and $m_2\in 2\N_0$ in the boundary-transmission conditions (\ref{normalm1})-(\ref{obliquem2}) and assume that (\ref{normalcondition}) holds. Then both operators ${\mathcal P}$ and  $W$ are normally solvable operators.
\end{teo}

\begin{pro}
For $m_1\in 2\N_0+1$ and $m_2\in 2\N_0$ the lifted Fourier symbol reads
$$
\Phi_0=\rho^{1-m_1-m_2}
\left[
   \begin{array}{cc}
-b_2^+-i^{m_2}\check{b}_2^+(\xi\beta^{-1})^{m_2}
  &b_2^-+i^{m_2}\check{b}_2^-(\xi\beta^{-1})^{m_2} \\  
-c_2^+-i^{m_2}\check{c}_2^+(\xi\beta^{-1})^{m_2} & 
c_2^-+i^{m_2}\check{c}_2^-(\xi\beta^{-1})^{m_2}
  \end{array}
 \right] , 
$$
which has no jumps at infinity, since both $\rho(\xi)^{1-m_1-m_2}$ and $(\xi\beta(\xi)^{-1})^{m_2}$ tend to one as $\xi$ tends to $\pm \infty$. Therefore, the corresponding WHO has always a closed image and, by the equivalence relation (\ref{equivalentrelation}), so does the operator ${\mathcal P}$.
\end{pro}

In the present paper we were able to achieve the image normalization of particular WHOs which arise from relevant boundary-transmission value problems for a junction of two half-planes. For theoretical and practical reasons it is most important to be able to answer further questions about the invertibility or the Fredholm properties of these operators. We plan to do this in a future work.

\vspace{5mm}

{\bf Aknowledgements:} This work was partially supported by {\it Funda\c{c}\~ao para a Ci\^{e}ncia e a Tecnologia} through {\it Centro de An\'alise Funcional e Aplica\c{c}\~oes}
of Instituto Superior T\'ecnico, Technical University of Lisbon, Portugal. 

%%%%%%%%%%%%%%%%%%%%%%%%

%%%%%%%%%%%%%%%%%%%%%%%%%%


\begin{thebibliography}{9}

 \bibitem{AmCastro04} L.P. Castro, A. Moura Santos,  {\it An operator approach for an oblique derivative boundary-transmission problem}, {Math. Meth.  Appl. Sci.} {\bf 27}  (2004)
        1469-1491.
        
        \bibitem{Duduchava79} R. Duduchava, in {\it Integral Equations with Fixed Singularities}, Teubner Texte zur Mathematik, Teubner, Leipzig, 1979.

        
        \bibitem{Eskin81} G.I. $\grave {E}$skin, in {\it Boundary Value Problems for Elliptic Pseudodifferential Equations},  Amer. Math. Soc., Providence, 1981.

       \bibitem{Kravchenko85} V.G. Kravchenko, {\it On normalization of singular integral operators}, {Soviet Math. Dokl.} {\bf 32} (1985) 880-883.
        
          \bibitem{MeisterSpeck89} E. Meister, F.-O. Speck, {\it Modern Wiener-Hopf methods in diffraction theory}, {Pitman Res. Notes Math. Ser.} {\bf 216}  (1989)
        130-171.

   \bibitem{MikhlinProssdorf86}   S.G. Mikhlin, S. Pr\"{o}ssdorf,    in {\it Singular Integral Operators}, Springer, Berlin, 1986.   
   
       \bibitem{AmSpeck97} A. Moura Santos, F.-O. Speck, {\it Sommerfeld diffraction problems with oblique derivatives}, {Math. Meth.  Appl. Sci.} {\bf 20} (1997)
        635-652. 
        
\bibitem{AmSpeckTeixeira98} A. Moura Santos, F.-O. Speck,  F.S. Teixeira, {\it Minimal normalization of Wiener-Hopf operators in spaces of Bessel potentials}, {J. Math. Anal. Appl.} {\bf 225}  (1998)
        501-531.
        
  \bibitem{Norris95} A.N. Norris, G.R. Wickham, {\it Accoustic diffraction from the junction of two flat plates}, {Proc. R. Soc. Lond. Ser. A} {\bf 451}  (1995)
        631-655.
        
        \bibitem{Rojas89} R.J. Rojas, P.H. Pathak,  {\it Diffraction of EM waves by a dielectric/ferrite half-plane and related configurations}, {IEEE Trans. Antennas and Propagation} {\bf 27}  (1989)
        751-763.
        
        \bibitem{SantosTeixeira95} P.A.  Santos,  F.S. Teixeira, {\it Sommerfeld half-plane problems with higher-order boundary  conditions}, {Math. Nachr.} {\bf 171} (1995)
        269-282.
        
        \bibitem{Serbest95} A.H. Serbest, E. L\"{u}neburg, {\it Diffraction at the junction of a two-impedance half-plane and a resistive half-plane}, {Antennas and Propagation Society International Symposium} {\bf 2} (1995)
        1356-1359.     
                    
            \bibitem{Speck86}  F.-O. Speck, {\it Mixed boundary value problems of the type of Sommerfeld's half-plane problem}, {Proc. Royal Soc. Edinburgh} {\bf 104A} (1986)
        261-277.


        \end{thebibliography}
\end{document}